\documentclass[a4,12pt,multicol,makeidx]{article}
\def\origin{%p'=(p_1,...,p_m), 
\vskip-\baselineskip\vskip-\topskip%
  \vbox to 0pt{\vskip-1in%
    \hbox to 0pt{\hskip-1in%
      \hbox to 0pt{\vrule width 1cm height .4pt depth 0mm\hss}%
      \vbox to 0pt{\hrule width .4pt height 0pt depth 1cm\vss}%
    \hss}%
  \vss}%%
  \vskip-\baselineskip
  \vbox to 0pt{\vskip-1in\vskip3cm%
    \hbox to 0pt{\hskip-1in\hskip3cm%
      \hbox to 0pt{\hss\vrule width 2cm height .4pt depth 0mm\hss}%
      \vbox to 0pt{\vss\hrule width .4pt height 1cm depth 1cm\vss}%
    \hss}%
  \vss}%
\vskip5mm\hskip10mm (3cm,3cm)
}%
\def\CA{\mathcal A}   
\def\a{\alpha} 
\def\d{\delta} \def\l{\lambda}   \def\p{\partial}  \def\e{\varepsilon} 
\def\n{\nabla}    
\def\leq{\underline{<}} 
\def\la{\langle} \def\ra{\rangle} \def\hx{\hat{x}} \def\hy{\hat{y}}
\def\ox{\overline{x}} \def\oy{\overline{y}}
\newenvironment{theorem}{%
\par \bigskip \it}{%
\bigskip \par}

\newenvironment{lemma}{%
\par \bigskip \it}{%
\bigskip \par}

\newenvironment{definition}{%
\par \bigskip \it}{%
\bigskip \par}

\newpage
\makeindex
\title{A localization of the   L{\'e}vy operators
 arising in mathematical finances.
}
\author{Mariko Arisawa\\ GSIS, Tohoku University
\\ Aramaki 09, Aoba-ku, Sendai 980-8579, JAPAN\\
E-mail: arisawa@math.is.tohoku.ac.jp
}
\date{}
\begin{document}
\maketitle
\bigskip
%\pagestyle{plain}

%%%%%Section 1.%%%%%
\newpage

\section{Introduction} 

$\quad$ We study the uniform H{\"o}lder continuity of the solutions of the following problem. 
$$
	F(x,\n v(x),\n^2 v(x))-\int_{{\bf R^N}}  
	[v(x+z)-v(x)
		\qquad\qquad\qquad
$$
\begin{equation}\label{first}
	-{\bf 1}_{|z|\leq 1}\la \n v(x),z \ra]c(z) dz -g(x) =0 \qquad x\in {\bf R^N},
\end{equation}
where 
 $c(z)dz$ is a positive Radon measure, called L{\'e}vy density, defined on
  ${\bf R^N}$ such that 
\begin{equation}\label{cintegral}
	\int_{{\bf R^N}} \min(|z|^2,1) c(z)dz \leq C_1 ,
\end{equation}
\begin{equation}\label{ccont}
	\frac{C_2}{|z|^{N+\gamma}}\leq |c(z)|\leq 
	\frac{C_3}{|z|^{N+\gamma}}\quad
	\forall z\in {\bf R^N}\cap \{|z|\leq 1\},
\end{equation}
where $\gamma\in (0,2)$, $C_i>0$ ($1\leq i\leq 3$) are constants.
 We assume that there exists a "uniform" constant $M>1$ such that
  for a constant $\theta_0\in [0,1]$, 
\begin{equation}\label{gcont}
	|g(x)-g(y)|\leq M|x-y|^{\theta_0} \quad \forall x,y\in {\bf R^N},
\end{equation}
and 
\begin{equation}\label{vbound}
	\sup_{x\in {\bf R^N}}|v|<M.
\end{equation}
The second-order fully nonlinear partial differential  operator $F$ is continuous in ${\bf R^N}\times$${\bf R^N}\times$${\bf S^N}$, and assumed to satisfy the following two conditions. \\
(Degenerate ellipticity) : 
$$
	F(x,p,X)\geq F(x,p,Y) \quad \hbox{if}\quad X\leq Y,
	\qquad\qquad
$$
\begin{equation}\label{delliptic}
	\qquad\qquad\qquad\quad\qquad\qquad
	\quad \forall x\in {\bf R^N}, \quad \forall p\in {\bf R^N},
	\quad \forall X,Y\in {\bf S^N}.
\end{equation}
(Continuity I) : 
There are modulus of continuity functions $w$ and 
$\eta$ from ${\bf R^+}\cup \{0\}\to {\bf R^+}\cup \{0\}$ such that 
$\lim_{\sigma\downarrow 0}w(\sigma)=0$, $\lim_{\sigma\downarrow 0}\eta(\sigma)=0$, and 
\begin{equation}\label{Fcont}
	|F(x,p,X)-F(y,p,X)|\leq w(|x-y|)|p|^q +\eta(|x-y|)||X||
	\qquad
\end{equation}
$$
	\qquad
	\qquad \forall x,y\in {\bf R^N},\quad \forall p\in {\bf R^N}, \quad \forall 	X\in {\bf S^N},
$$	
where $q\geq 1$. \\

	We study this problem in the framework of the viscosity solutions for the integro-differential equations, the definition of which is introduced in Arisawa \cite{ar4} (see also \cite{ar5} and \cite{ar6}). 
	The definition is the following. 
 In  order to get rid of the singularity of the L{\'e}vy measure, we  shall use
 the following superjet (resp. subjet) and its residue. Let 
$\hx \in {\bf R^N}$, and let $(p,X)\in J^{2,+}_{{\bf R^N}}u(\hx)$ (resp. $(p,X)\in J^{2,-}_{{\bf R^N}}u(\hx)$) be a second-order superjet (resp. subjet) of $u$ at $\hx$. Then, for any $\d>0$ there exists $\e>0$ such that 
\begin{equation}\label{edp}	
	u(\hx+z)\leq u(\hat{x})+\la p,z \ra 
	+\frac{1}{2}\la X z,z \ra +\delta |z|^2  \qquad \hbox{if}\quad |z|\leq \e
\end{equation}
(resp.
\begin{equation}\label{edm}	
		v(\hx+z)\geq v(\hat{x})+\la p,z \ra 
	+\frac{1}{2}\la X z , z \ra 
	-\delta |z|^2  \qquad \hbox{if}\quad |z|\leq \e
\end{equation}
 ) holds. 
 We use this pair of numbers $(\e,\d)$ satisfying 
(\ref{edp}) (resp. (\ref{edm})) for any $(p,X)\in J_{{\bf R^N}}^{2,+}u(\hx)$ 
(resp. $(p,X)\in J_{{\bf R^N}}^{2,-}v(\hx)$) in the following definition of 
viscosity solutions.\\
	\begin{definition}{\bf Definition 1.1.}  Let $u\in USC({\bf R^N})$ (resp. $v\in LSC({\bf R^N})$). We say that $u$ (resp. $v$) is a viscosity subsolution (resp. supersolution) of (\ref{first}), if for any $\hx\in {\bf R^N}$, any $(p,X)\in J_{{\bf R^N}}^{2,+}u(\hx)$ (resp. $\in J_{{\bf R^N}}^{2,-}v(\hx)$), and 
any pair of numbers $(\e,\delta)$ satisfying  (\ref{edp}) (resp.(\ref{edm})), 
the following holds for any $0< \e'\leq \e$
$$
	F(\hat{x},p,X)
	-\int_{|z|<\e'} \frac{1}{2}\la(X+2\d I)z,z \ra c(z)dz
$$
$$
	- \int_{|z|\geq \e'} [u(\hat{x}+z)-u(\hat{x})
	-{\bf 1}_{|z|\leq 1}\la z,p\ra] c(z)dz \leq 0.
$$
(resp.
$$
	F(\hat{x},p,X)
	-\int_{|z|<\e'} \frac{1}{2}\la (X-2\d I)z,z \ra c(z)dz
$$
$$
	-\int_{|z|\geq \e'} [v(\hat{x}+z)-v(\hat{x})
	-{\bf 1}_{|z|\leq 1}\la z,p\ra] c(z)dz \geq 0.
$$
If $u$ is both a viscosity subsolution and a viscosity supersolution , it is called a viscosity solution.
\end{definition}

In the framework of the viscosity solutions in Definition 1.1, we have the existence and the comparison results in \cite{ar4}, \cite{ar5} and \cite{ar6}. For the convenience of the readres, we shall give typical comparison results and the proof in \S2 in below.  \\

 Then, we claim the uniform H{\"o}lder continuity of $u$ in the following two cases.\\
(I) $N=1$. \\
(II) $N\geq 2$, and $F$ satisfies the following  uniform ellipticity. \\
(Uniform ellipticity) : There exists $\lambda_0>0$ such that 
$$
	F(x,p,X)-F(x,p,Y)\geq \lambda_0 (Y-X) \quad \hbox{if}\quad
	X\leq Y,
$$
\begin{equation}\label{uelliptic}
	\qquad\qquad\qquad\qquad\qquad\qquad\qquad
	\forall x\in {\bf R^N},
	\quad \forall p\in {\bf R^N},\quad \forall X,Y\in{\bf S^N}.
\end{equation}  
In the case of (I), we claim that for any 
  $\theta\in (0,\min\{1,\theta_0+\gamma\})$ there  exists $C_{\theta}>0$ such that
\begin{equation}\label{vhold}
	|v(x)-v(y)|\leq C_{\theta}|x-y|^{\theta}\quad \forall x,y\in {\bf R^N},
\end{equation}
where $C_{\theta}>0$ depends only on $M$ and $C_1$. (See Theorem 3.1 in below.) In the case of (II), we claim that for any 
$\theta\in (0,1)$, there exists $C_{\theta}>0$ such that (\ref{vhold}) holds. (See Theorem 3.2 in below.) (These results hold for more general problem 
$$
	F(x,\n v(x),\n^2 v(x))+\sup_{\a\in \CA}\{-\int_{{\bf R^N}}  
	[v(x+z)-v(x)
		\qquad\qquad\qquad
$$
$$
	-{\bf 1}_{|z|\leq 1}\la \n v(x),z \ra]c(x,z,\a) dz -g(x,\a)\} =0 \qquad x\in {\bf R^N},
$$
which we do not treat here.)\\

	As for the case  other than (I) and (II), that is $N\geq 2$ and  $F$ is not necessarily uniformly elliptic (i.e. (\ref{uelliptic}) is not satisfied),  
we study the following two problems in the torus ${\bf T^N}$ instead of (\ref{first}). The first one is, for $\l>0$,   
$$
	\l v(x)+ H(\n v(x))-\int_{{\bf R^N}}  
	[v(x+z)-v(x)
	\qquad\qquad\qquad\qquad\qquad
$$
\begin{equation}\label{lambH}
	-{\bf 1}_{|z|\leq 1}\la \n v(x),z \ra]c(z) dz -g(x)=0 \qquad x\in {\bf T^N}.
\end{equation}
And the second one is  
$$
	\l v(x)+ F(x,\n v(x),\n^2 v(x)) -\int_{{\bf R^N}}  
	[v(x+z)-v(x)
	\qquad\qquad\qquad
$$
\begin{equation}\label{lamb}
	-{\bf 1}_{|z|\leq 1}\la \n v(x),z \ra]c(z) dz -g(x)=0 \qquad x\in {\bf T^N},
\end{equation}
where $\l>0$. Here $H$ is a first-order nonlinear operator, and $F$ is a fully 
nonlinear degenerate elliptic operator, satisfying the following conditions.\\
(Periodicity) : 
$$
	H(\cdot,p),\quad F(\cdot,p,X),\quad \hbox{and}\quad g(\cdot)\quad
	\hbox{are periodic in}\quad x\in {\bf T^N},\quad 
$$
\begin{equation}\label{periodicity}	
	\qquad\qquad\qquad\qquad\qquad\qquad
	\hbox{for}\quad
	\forall p\in {\bf R^N},
	\quad \forall X\in {\bf S^N}.
\end{equation}

\leftline{(Partial uniform ellipticity) : 
There exists a constant $\l_1> 0$ such that }
\begin{eqnarray}\label{A1}
	&F(x,p,X) \geq F(x,p,Y) +\l_1 Tr (Y'-X') 
	\quad \forall x\in {\bf T^N}, 
	\quad \forall p\in {\bf R}^N, \nonumber\\
	&\quad 
	\forall X,Y\in {\bf S}^N,
	 \quad X= \left(
		\begin{array}{cc}
		X'& X_{12}\\ 
		X_{21} & X_{22}
		\end{array}
		\right), \quad Y= \left(
		\begin{array}{cc}
		Y'& Y_{12}\\
		Y_{21} & Y_{22}
		\end{array}
		\right),
\end{eqnarray}
$$
	\hbox{where} \quad X'\leq Y' 
	(X',Y'\in {\bf S}^M),\quad  0<M \leq N.
	\qquad\qquad\qquad\qquad\qquad\qquad\qquad\qquad\qquad\qquad\quad
$$

(Continuity II) : There are modulus of continuity functions $w'$ and 
$\eta'$ from ${\bf R^+}\cup \{0\}\to {\bf R^+}\cup \{0\}$ such that 
$\lim_{\sigma\downarrow 0}w'(\sigma)=0$, $\lim_{\sigma\downarrow 0}\eta'(\sigma)=0$, and 
\begin{eqnarray}
	|F(x,p,X)-F(y,p,X)| \leq w'(|x-y|)|p'|^{q'} +\eta'(|x-y|)||X'|| \qquad
	\nonumber \\
	\forall x,y \in {\bf T^N}, 
	\quad \forall p=(p',p'')\in {\bf R}^M\times {\bf R}^{m},\quad
	\forall X= \left(
		\begin{array}{cc}
		X'& X_{12}\\
		X_{21} & X_{22}
		\end{array}
		\right) \in {\bf S}^N, \nonumber 
\end{eqnarray}
\begin{equation}\label{pFcont}
	\hbox{where}\quad X'\in {\bf S}^M,\quad M+m=N, \quad q'\geq 1. 
	\qquad\qquad\qquad\qquad\qquad\qquad
\end{equation}
 Roughly speaking, we claim that for any 
 $\theta\in (0,\theta_0)$ ($\theta_0>0$), 
there exists $C_{\theta}>0$  such that 
\begin{equation}\label{lambvhold}
	|v(x)-v(y)|\leq \frac{C_{\theta}}{\lambda}|x-y|^{\theta}\quad \forall x,y\in {\bf T^N},
\end{equation}
where $C_{\theta}>0$ is independent on $\l>0$. (See Theorems 4.1 and 4.2 in below.)
	$\quad$The method to derive the above uniform H{\"o}lder continuity (\ref{vhold}) and the H{\"o}lder continuity (\ref{lambvhold}) is based on the argument 
used in the proof of  the comparison result. (See Ishii and Lions \cite{ipl}, for the similar argument in the PDE case.) \\
	
	Next, we shall state the strong maximum principle for the L{\'e}vy operator. In \cite{menal}, 
 for the second-order uniformly elliptic integro-differential 
 operator
$$
 	-\sum_{i,j=1}^N \overline{a}_{ij}\frac{\p^2 v}{\p x_i \p x_j}-\sum_{i=1}^N
 	\overline{b}_i\frac{\p v}{\p x_i}-\int_{{\bf R^N}}[v(x+z)
	-v(x)-\la \n v(x),z \ra]c(x,z) dz  
$$
\begin{equation}\label{ergod}
	\qquad\qquad\qquad\qquad\qquad\qquad\qquad\qquad\qquad\qquad\qquad
	\qquad x\in {\bf R^N},
\end{equation}
 the strong maximum principle was given, 
 where $\lambda_0 I \leq (\overline{a}_{ij})_{1\leq i,j\leq N} \leq \Lambda_0 I$ ($0<\lambda_0\leq \Lambda_0$). See also, Cancelier \cite{can} for another type of the maximum principle. Here, we shall give the strong maximum principle in ${\bf R^N}$ without assuming the uniform ellipticity of the partial differential operator $F$ in (\ref{first})  (see Theorem 5.1 in below, and M. Arisawa and P.-L. Lions \cite{mpl2}). \\
	
	Finally, we shall apply these regularity results (\ref{vhold}), (\ref{lambvhold}) and the strong maximum principle, to study the so-called  ergodic problem. In the case of the Hamilton-Jacobi-Bellman (HJB) operator
$$
	\sup_{\a\in \CA}\{-\sum_{i,j=1}^{N} a_{ij}(x,\a)\frac{\p^2 u}{\p x_i\p x_j}
	-\sum_{i=1}^{N} b_i(x,\a)\frac{\p u}{\p x_i}-f(x,\a)
	\},
$$
the ergodicity of the corresponding controlled diffusion process, for example in the torus ${\bf T^N}={\bf R^N}\backslash {\bf Z^N}$, can be studied by the existence of a unique real number $d_f$ such that the following problem admits a periodic viscosity 
solution $u$ : 
$$
	d_f+\sup_{\a\in \CA}\{-\sum_{i,j=1}^{N} a_{ij}(x,\a)\frac{\p^2 u}{\p x_i\p x_j}
	-\sum_{i=1}^{N} b_i(x,\a)\frac{\p u}{\p x_i}-f(x,\a)
	\}=0\quad x\in {\bf T^N}.
$$
We refer the readers to M. Arisawa and P.-L. Lions \cite{apl}, M. Arisawa \cite{ar1}, \cite{ar2}, 
for more details. From the analogy of the diffusion case, 
here we shall formulate the ergodic problem for the integro-differential equations as follows. \\

(Ergodic problem) Is there a unique number $d_f$ depending only on $f(x)$ such that the following problem has a periodic viscosity solution $u(x)$ defined on $\bf T^N$ ? 
$$
	d_f+F(x,\n u,\n^2 u) -\int_{{\bf R^N}}  
	[u(x+z)
	\qquad\qquad\qquad\qquad\qquad\qquad
$$
$$
	-u(x)-{\bf 1}_{|z|\leq 1}\la \n u(x),z \ra] c(z)dz -f(x)=0 \qquad x\in {\bf T^N}.
$$
 The results on the existence of the above number $d_f$ is stated in Theorem 6.1 in below. \\
%%%%%%%%%%%%%%%%%%%%%%%%%%%%%%%%%%%%%%%%%%%%%%%%%%%%%%%%%%%%%%%%%%%%%%%%%%%%%
\section{Comparison results} In this section, we give some typical comparison results for the integro-differential equations in the framework of the solution in Definition 1.1. We consider 
$$
	\l u + F(x,\nabla u,\nabla^2 u) -\int_{{\bf R^N}}
	u(x+z)-u(x)\qquad \qquad
$$
\begin{equation}\label{problemI}
	\qquad -{\bf 1}_{|z|\leq 1}\la z,\n u(x) \ra q(dz)=0 \qquad \hbox{in} \quad \Omega,
\end{equation}
where $\l>0$, and $\Omega$ is a bounded domain in ${\bf R^N}$,  
with either the Dirichlet B.C.:
\begin{equation}\label{dirichletI}
	u(x)=g(x) \quad \forall x\in \Omega^c,
\end{equation}
or the Periodic B.C.:
\begin{equation}\label{periodic}
	\Omega={\bf T^N}={\bf R^N}\backslash {\bf Z^N},\quad
	u(x)\quad\hbox{is periodic in}\quad {\bf T^N},
\end{equation}
where $g$ is a given continuous function in ${\Omega^c}$. The second-order partial defferential operator $F$ is degenerate elliptic, which satisfies

(Degenerate ellipticity) (cf. \cite{users} (3.14)): 
There exists a function $w(\cdot)$$:[0,\infty)\to[0,\infty)$, $w(0+)=0$ 
such that
\begin{equation}\label{a2}
	F(y,r,p,Y)-F(x,r,p,X)\leq w(\alpha|x-y|^2+|x-y|(|p|+1))
\end{equation}
$$
	\qquad \qquad \qquad\qquad\qquad
	\hbox{for}
	\quad x,y\in\overline{\Omega}, \quad r\in {\bf R}, \quad p\in {\bf R^N}
$$
for any $\a>0$, and for any $X,Y\in {\bf S^N}$ such that 
\begin{equation}\label{three}
	-3\alpha \left(
		\begin{array}{cc}
			I&O\\
			O&I
		\end{array}
	\right) \leq 
		\left(
		\begin{array}{cc}
			X&O\\
			O&-Y
		\end{array}
	\right) \leq 
	3\alpha \left(
		\begin{array}{cc}
			I&-I\\
			-I&I
		\end{array}
	\right).
\end{equation}

 For the above, we have the following results. 

{\bf Theorem 2.1}\begin{theorem} Assume that $\Omega$ is bounded, and that (\ref{cintegral}), (\ref{a2}) hold. Let $u\in USC({\bf R^N})$ and $v\in LSC(\bf R^N)$ be respectively a viscosity subsolution and a supersolution of (\ref{problemI}) in $\Omega$, which satisfy $u\leq v$ on $\Omega^c$. Then, 
$$
	u\leq v \qquad \hbox{in}\quad \Omega.
$$
\end{theorem}
 
{\bf Theorem 2.2}\begin{theorem} Let $\Omega={\bf T^N}$. Assume that (\ref{cintegral}), (\ref{a2}) hold and that $F$ is periodic in $x\in {\bf T^N}$. Let $u\in USC({\bf T^N})$ and $v\in LSC(\bf T^N)$ be respectively a viscosity subsolution and a supersolution of (\ref{problemI}) in $\Omega$. Then, 
$$
	u\leq v \qquad \hbox{in}\quad \Omega.
$$
\end{theorem}

{\bf Remark 2.1} The above comparison results hold in more general situations. 
For example, $\Omega$ can be ${\bf R^N}$ by assuming that $u$ and $v$ are bounded, or the nonlocal operator can be in the form of 
$$
	-\int_{\{z\in {\bf R^N}|x+z\in \overline{\Omega}\}}
	[u(x+z)-u(x)-{\bf 1}_{|z|\leq 1}\la z,\n u(x) \ra ]c(z)dz, 
	\quad x\in {\bf \Omega},
$$
with the Neumann type boundary condition on $\p\Omega$, etc... We refer the readers to 
\cite{ar4}, \cite{ar5} and \cite{ar6}. \\

	In order to prove the above claims, we use the following two Lemmas. (See \cite{ar6}.) 
The first Lemma is the approximation by the supconvolution and the infconvolution.\\
	
{\bf Lemma 2.3}\begin{theorem}  Let $u$ and $v$ be respectively a bounded viscosity subsolution and a bounded supersolution of (\ref{problemI}). Define for $r>0$, the supconvolution $u^r$ and the infconvolution $v_r$ of $u$ and $v$ as follows.
\begin{equation}\label{supconvolution}
	u^r(x)=\sup_{y\in {\bf R^N}}\{ 
	u(y)-\frac{1}{2r^2}|x-y|^2
	\}\qquad (\hbox{supconvolution}).
\end{equation}
\begin{equation}\label{infconvolution}
	v_r(x)=\inf_{y\in {\bf R^N}}\{ 
	v(y)+\frac{1}{2r^2}|x-y|^2
	\}\qquad (\hbox{infconvolution}).
\end{equation}
Then, for any $\nu>0$ there exists $r>0$ such that $u^r$ and $v_r$ are respectively a subsolution and a supersolution of the following problems. 
$$
	\l u_r + F(x,\nabla u^r,\nabla^2 u^r) -\int_{{\bf R^N}}
	u^r(x+z)-u^r(x)\qquad \qquad
$$
\begin{equation}\label{subproblem}
	\qquad -{\bf 1}_{|z|\leq 1}\la z,\n u^r(x) \ra q(dz)\}\leq \nu \qquad \hbox{in} \quad \Omega_r.
\end{equation}
$$
	\l v_r+ F(x,\nabla v_r,\nabla^2 v_r) -\int_{{\bf R^N}}
	v_r(x+z)-v_r(x)\qquad \qquad
$$
\begin{equation}\label{superproblem}
	\qquad -{\bf 1}_{|z|\leq 1}\la z,\n v_r(x) \ra q(dz)\}\geq -\nu \qquad \hbox{in} \quad \Omega_r,
\end{equation}
where $\Omega_r=\{x\in {\Omega}|\quad dist(x,\p\Omega)> \sqrt{2M}r\}$ for $M=\max\{\sup_{\overline{\Omega}}|u|,\sup_{\overline{\Omega}}|v|\}$. 
\end{theorem}

Remark that $u^r$ is semiconvex, $v_r$ is semiconcave, and both are Lipschitz continuous in ${\bf R^N}$. The second lemma comes 
from the Jensen's maximum principle and the Alexandrov's theorem (see \cite{users} and \cite{fleming}). The last claim of this lemma is quite important in the limit procedure in the nonlocal term. 
\begin{theorem}{\bf Lemma 2.4}
	Let $U$ be semiconvex and $V$ be semiconcave in $\Omega$. 
	For 
$\phi(x,y)=\a |x-y|^2$ ($\a>0$) consider 
$\Phi(x,y)=U(x)-V(y)-\phi(x,y)$, 
and assume that $(\ox,\oy)$ is an interior maximum of $\Phi$ 
in $\overline{\Omega}\times \overline{\Omega}$. Assume also that there is an open precompact subset $O$ of $\Omega\times \Omega$ containing $(\ox,\oy)$,  and that 
$\mu$
$=\sup_{O} \Phi(x,y) -\sup_{\p O} \Phi(x,y) >0$.
Then, the following holds.\\
(i) There exists a sequence of points $(x_m,y_m)\in O$ ($m\in {\bf N}$) 
such that $\lim_{m\to \infty} (x_m,y_m)=(\ox,\oy)$, and 
$(p_m,X_m)\in J^{2,+}_{\Omega}U(x_m)$,  $(p_m',Y_m)\in J^{2,-}_{\Omega}V(y_m)$
 such that 
$\lim_{m\to \infty}p_m$$=\lim_{m\to \infty}p'_m$$=2\a(x_m-y_m)=p$, and 
$X_m\leq Y_m\quad \forall m$.\\
(ii) For $P_m=(p_m-p,-({p'}_m-p))$, $\Phi_m(x,y)=\Phi(x,y)-\la P_m,(x,y)\ra$
 takes a maximum at $(x_m,y_m)$ in O. \\
(iii) The following holds for any $z\in {\bf R^N}$ such that
 $(x_m+z,y_m+z)\in O$.
\begin{equation}\label{important}
	U(x_m+z)-U(_m)-\la p_m,z\ra
	\leq V(y_m+z)-V(y_m)-\la p'_m,z\ra. 
\end{equation}
\end{theorem}

	We admit the above claims here. (In fact, the proofs of Lemma 2.3 and 2.4 are  not so difficult, see for example \cite{users} and \cite{fleming}.)\\

$Proof\quad of\quad Theorem\quad 2.1.$ 
We use the argument by contradiction, and assume that $\max_{\overline{\Omega}}(u-v)$$=(u-v)(x_0)=M_0>0$ for $x_0\in \Omega$. Then, we approximate $u$ by $u^r$ (supconvolution) and $v$ by $v_r$ (infconvolution), which are a subsolution  and a supersolution of (\ref{subproblem}) and (\ref{superproblem}), respectively. Clearly, $\max_{\overline{\Omega}}(u^r-v_r)$$\geq M_0>0$. Let $\ox\in \Omega$ be the maximizer of $u^r-v_r$. In the following, we abbreviate the index and  write $u=u^r$, $v=v_r$ without any confusion. As in the PDE theory, consider $\Phi(x,y)=u(x)-v(y)-\a|x-y|^2$, and let $(\hx,\hy)$  be the maximizer of $\Phi$. Then, from Lemma 2.3 there exists $(x_m,y_m)\in \Omega$ ($m\in {\bf N}$) 
such that $\lim_{m\to \infty} (x_m,y_m)=(\hx,\hy)$, and we can take 
$(\e_m,\d_m)$  a pair of positive numbers such that
$u(x_m+z)\leq u(x_m)+\la p_m,z\ra+\frac{1}{2}\la X_m z,z\ra+\d_m|z|^2$, 
$v(y_m+z)\geq v(y_m)+\la p'_m,z\ra+\frac{1}{2}\la Y_m z,z\ra-\d_m|z|^2$, 
for $\forall |z|\leq \e_m$. 
	From the definition of the viscosity solutions, we have 
$$
	F(x_m,u(x_m),p_m,X_m) -\int_{|z|\leq \e_m}
	\frac{1}{2}\la (X_m+2\d_m I)z,z \ra dq(z)
$$
$$
	-\int_{|z|\geq \e_m}
	u(x_m+z)-u(x_m)
	-{\bf 1}_{|z|\leq 1}\la z,p_m \ra q(dz)\leq \nu,
$$
$$
	F(y_m,v(y_m),p'_m,Y_m) -\int_{|z|\leq \e_m}
	\frac{1}{2}\la (Y_m-2\d_m I)z,z \ra dq(z)
$$
$$
	-\int_{|z|\geq \e_m}
	v(y_m+z)-v(y_m)
	-{\bf 1}_{|z|\leq 1}\la z,p'_m \ra q(dz)\geq -\nu.
$$	
By taking the difference of the above two inequalities,  by using (\ref{important}), and by passing $m\to \infty$ (thanking to (\ref{important}), it is available), we can obtain the desired contradiction. The claim 
$u\leq v$ is proved. \\

{\bf Remark 2.2} As for the usage of (\ref{important}) in the  limit procedure 
$m\to \infty$ in the proof of Theorem 2.1, we refer the interested readers 
 to the similar argument in the  proof of Theorem 3.2 in below. \\

%%%%%%%%%%%%%%%%%%%%%%%%%%%%%%%%%%%%%%%%%%%%%%%%%%%%%%%%%%%%%%%%%%%%%%%%%%%%%%
\section{Uniform H{\"o}lder continuities of viscosity solutions}

$\quad$ In this section, we study the uniform H{\"o}lder continuities of viscosity solutions of (\ref{first}) in the cases of (I) and (II). 
\\

{\bf Theorem 3.1.$\quad$}
\begin{theorem} Let $N=1$, and let $v$ be a viscosity solution of (\ref{first}) satisfying (\ref{vbound}). Assume that (\ref{cintegral}), (\ref{ccont})  and (\ref{gcont}) 
hold, where $\gamma\in (0,2)$. 
Assume also that  $F$ satisfies (\ref{delliptic}) and (\ref{Fcont}), where 
 there exist constants $L>0$, $\rho_i>0$ ($i=1,2$) such that
\begin{equation}\label{os211}
	\lim_{s\downarrow 0}w(s)s^{-\rho_1}\leq L,\quad 
	\lim_{s\downarrow 0}\eta(s)s^{-\rho_2}\leq L,
\end{equation}
and $\rho_1+\gamma>q$, $\rho_2+\gamma>2$. 
Then for any $\theta\in (0,\min\{1,\theta_0+\gamma\})$, there exists a constant $C_{\theta}>0$ such that (\ref{vhold}) holds. The constant $C_{\theta}$ 
depends only on $M>0$ and $C_i$ ($1\leq i\leq 3$).
\end{theorem}

{\bf Theorem 3.2.$\quad$}
\begin{theorem} Let $N\geq 2$, and let $v$ be a viscosity solution of (\ref{first}) satisfying (\ref{vbound}). Assume that (\ref{cintegral}), (\ref{ccont}) 
and (\ref{gcont})
hold, where $\gamma\in (0,2)$.

If $F$ satisfies (\ref{Fcont}) and (\ref{uelliptic}), 
where there exist constants $L>0$, $\rho_1>0$ such that 
\begin{equation}\label{os212}
	\lim_{s\downarrow 0}w(s)s^{-\rho_1}\leq L,
\end{equation}
and $\rho_1 + 2>q$, 
then for any $\theta\in (0,1)$, there exists a constant $C_{\theta}>0$ such that (\ref{vhold}) holds. 
The constant $C_{\theta}$ 
depends only on $M>0$ and $C_i$ ($1\leq i\leq 3$).
\end{theorem}

	The following lemma gives the relationship between $\d$ and $\e$ in Definition 1.1, and is used in the proofs in below. \\

{\bf Lemma 3.3.\quad} 
\begin{theorem} Let $\phi(z)=C_{\theta}|z|^{\theta}$ ($z\in {\bf R^N}$), $\theta\in (0,1)$, $r>0$, and 
let $\hat{z}\in \{z\in {\bf R^N}|\quad |z|< r,\quad z\neq 0\}$ be fixed. Then, there exists $\overline{C}>0$ 
such that for any $\d>0$, and for any $z\in {\bf R^N}$ such that $|z|\leq \frac{|\hat{z}|}{2}$, if $z$ satisfies 
\begin{equation}\label{z}
	|z|\leq \d \overline{C} |\hat{z}|^{3-\theta},
\end{equation}
we have 
\begin{equation}\label{ldevelop}
	|\phi(\hat{z}+z)-\phi(\hat{z})-\la \n \phi(\hat{z}), z\ra 
	-\frac{1}{2}\la \n^2 \phi(\hat{z})z,z\ra| \leq \d|z|^2. 
\end{equation}
The constant $\overline{C}$ is independent on $r$, $\theta$, and $\hat{z}$. 
\end{theorem}

$Proof$ $of$ $Lemma$ $3.3.\quad$ From the Taylor expansion of $\phi$ at $\hat{z}$
$$
	\phi(\hat{z}+z)-\phi(\hat{z})-\la \n \phi(\hat{z}), z\ra 
	-\frac{1}{2}\la \n^2 \phi(\hat{z})z,z\ra 
	=\frac{1}{3!}\sum_{i,j,k=1}^N \frac{\p^3\phi(\hat{z}+\rho(z)z)}{\p z_i \p z_j \p z_k}
	z_iz_jz_k 
$$
$$
	\qquad\qquad\qquad\qquad\qquad \hbox{for}\quad 
	z\in \{z\in {\bf R^N}|\quad|z|<\frac{|\hat{z}|}{2}\},
$$
where $\rho=\rho(z)\in (0,1)$. By calculating $\frac{\p^3 \phi}{\p x_i \p x_j \p x_k}$, we see that there exists a constant $C>0$ independent on $r$, $\theta$, $\d$ such that 
$$
	|\phi(\hat{z}+z)-\phi(\hat{z})-\la \n \phi(\hat{z}), z\ra 
	-\frac{1}{2}\la \n^2 \phi(\hat{z})z,z\ra|
	\leq C|\hat{z}+\rho z|^{\theta-3}|z|^3
$$
$$
	\qquad\qquad\qquad\qquad\qquad \hbox{for}\quad 
	z\in \{z\in {\bf R^N}|\quad|z|<\frac{|\hat{z}|}{2}\}.
$$
Then, if $|z|\leq \frac{\d}{C}|\hat{z}+\rho z|^{3-\theta}$
\begin{equation}\label{if}
	C|z|^3|\hat{z}+\rho z|^{\theta-3}\leq \d |z|^2.
\end{equation}
Since for $|z|\leq \frac{|\hat{z}|}{2}$, 
$$
	\frac{1}{2}|\hat{z}|\leq |\hat{z}+\rho z|\leq 2|\hat{z}|, 
$$
there exists $\overline{C}>0$ independent on $r$, $\theta$, $\d$, and $\rho$ such  that, if
$$
	\qquad\qquad|z|\leq \d \overline{C}|\hat{z}|^{3-\theta}\leq 
	\frac{\d}{C} |\hat{z}+\rho z|^{3-\theta},\qquad\qquad\qquad\qquad\qquad\qquad\qquad (\ref{z})'
$$
 then (\ref{if}) holds. 
Therefore, if $z$ satisfies (\ref{z}) with the above $\overline{C}>0$, and if 
$|z|<\frac{|\hat{z}|}{2}$, the inequality (\ref{ldevelop}) holds.\\
$$\quad$$

$Proof$ $of$ $Theorem$ $3.1.$  Fix an arbitrary number $\theta\in (0,1)$. Let
 $r_0>0$ be a small enough number which will be determined in the end of the proof. 
For  
$C_{\theta}>0$ such that 
\begin{equation}\label{ctheta}
	C_{\theta}r_0^{\theta}= 2M, 
\end{equation}
 we shall prove (\ref{vhold}), by the contradiction's argument. 
For $x,y\in {\bf R^N}$ such that $|x-y|\geq r_0$, from (\ref{vbound}) we have 
$$
	|v(x)-v(y)|\leq 2M\leq C_{\theta}|x-y|^{\theta}.
$$
Assume that there exist $x',y'\in {\bf R^N}$ ($|x'-y'|<r_0$) such that 
$$
	|v(x')-v(y')|>C_{\theta}|x'-y'|^{\theta},
$$
and we shall look for a contradiction.  
Consider for $\tau\in (0,1)$ 
$$
	\Phi(x,y)=v(x)-v(y)-C_{\theta}|x-y|^{\theta}-\frac{\tau}{2}|x|^2,
$$
and let $(\hat{x},\hat{y})$ be a maximum point of $\Phi$. Let us write $\phi(x,y)=C_{\theta}|x-y|^{\theta}$, and calculate 
$$
	\n_x \phi(x,y)=C_{\theta} \theta |x-y|^{\theta-2}(x-y)=-\n_y \phi(x,y)
$$
$$
	\n^2_{xx} \phi(x,y)=C_{\theta}\theta |x-y|^{\theta-2}I+
	C_{\theta}\theta(\theta-2) |x-y|^{\theta-4}(x-y)\otimes(x-y)
	=\n^2_{yy} \phi(x,y).
$$
Put $p=\n_x \phi(\hx,\hy)=-\n_y \phi(\hx,\hy)$, and $Q=\n^2_{xx} \phi(\hx,\hy)=\n^2_{yy} \phi(\hx,\hy)$. 
  Since 
$$
	\Phi(\hat{x}+z,\hat{y})=v(\hat{x}+z)-v(\hat{y})
	-C_{\theta}|\hat{x}+z-\hat{y}|^{\theta}-\frac{\tau}{2}(\hat{x}+z)^2
	\qquad\qquad\qquad
$$
$$
	\leq \Phi(\hat{x},\hat{y})=v(\hat{x})-v(\hat{y})
	-C_{\theta}|\hat{x}-\hat{y}|^{\theta}-\frac{\tau}{2}\hat{x}^2,
$$
for any $\d>0$ there exists $\e>0$ such that 
\begin{equation}\label{vdev1}
	v(\hat{x}+z)-v(\hat{x})\leq 
	C_{\theta}|\hat{x}+z-\hat{y}|^{\theta}-C_{\theta}|\hat{x}-\hat{y}|^{\theta}
	+\frac{\tau}{2}(\hat{x}+z)^2-\frac{\tau}{2}\hat{x}^2
\end{equation}
$$
	\qquad\qquad
	\leq ( p+\tau \hat{x})z  + \frac{1}{2}(Q+\tau)z^2 +\d z^2
	\quad \hbox{for}\quad |z|<\e. 
$$
Samely, since 
$$
	\Phi(\hat{x},\hat{y}+z)=v(\hat{x})-v(\hat{y}+z)
	-C_{\theta}|\hat{x}-(\hat{y}+z)|^{\theta}-\frac{\tau}{2}\hat{x}^2\qquad
	\qquad\qquad
$$
$$
	\qquad\qquad
	\leq \Phi(\hat{x},\hat{y})=v(\hat{x})-v(\hat{y})
	-C_{\theta}|\hat{x}-\hat{y}|^{\theta}-\frac{\tau}{2}\hat{x}^2,
$$
for any $\d>0$ there exists $\e>0$ such that 

\begin{equation}\label{vdev2}
	v(\hat{y}+z)-v(\hat{y})\geq 
	-(C_{\theta}|\hat{x}-(\hat{y}+z)|^{\theta}
	-C_{\theta}|\hat{x}-\hat{y}|^{\theta})
	\qquad\qquad\qquad
\end{equation}
$$
	\qquad\qquad
	\geq -(-pz  + \frac{1}{2}Qz^2 ) -\d z^2
	=pz  + \frac{1}{2}(-Q)z^2  -\d z^2
	\quad \hbox{for}\quad |z|<\e. 
$$
From the definition of viscosity solutions, 
by using the pair of numbers $(\e,\d)$ in (\ref{vdev1}) and (\ref{vdev2}), 
we have 
$$
	F(\hx,p+\tau \hx,Q+\tau)
 	-\int_{|z|\leq \e}\frac{1}{2}(Q+\tau+2\d)z^2 c(z)dz 	 \qquad\qquad
$$
$$
 	-\int_{|z|\geq \e} [v(\hx+z)-v(\hx)-{\bf 1_{|z|\leq 1}}(p+\tau \hx)z]
 	c(z)dz-g(\hx) \leq 0,
$$
 and 

 $$
 	F(\hy,p,-Q)-\int_{|z|\leq \e}\frac{1}{2}(-Q-2\d)z^2 c(z)dz	\qquad\qquad\qquad\qquad
$$
$$
	-\int_{|z|\geq \e} [v(\hy+z)-v(\hy)-{\bf 1_{|z|\leq 1}}pz ]
 	c(z)dz 
 	-g(\hy) \geq 0.
$$
 By taking the difference of the above two inequalities, 
 we have the following. 

$$
	F(\hx,p+\tau \hx,Q+\tau)-F(\hy,p,-Q)\qquad\qquad\qquad\qquad\qquad\qquad\qquad\qquad
$$
$$
	-\frac{1}{2}\int_{|z|\leq \e} (Q+\tau+2\d)z^2 
	c(z)dz
	-\frac{1}{2}\int_{|z|\leq \e} (Q+2\d)z^2  
	c(z)dz\qquad
$$
$$
	-\int_{|z|\geq \e} [v(\hx+z)-v(\hx)-{\bf 1_{|z|\leq 1}}(p+\tau \hx)z]
 	c(z)dz \qquad\qquad\qquad\qquad\qquad\qquad
$$
\begin{equation}\label{subst1}
	+\int_{|z|\geq \e} [v(\hy+z)-v(\hy)-{\bf 1_{|z|\leq 1}}pz]
 	c(z)dz 
 	\leq g(\hx)-g(\hy)+\nu.
\end{equation}
We need the estimates.  \\

{\bf Lemma 3.4.}
\begin{lemma}
The inequalities (\ref{vdev1}) and (\ref{vdev2}) hold with 
\begin{equation}\label{23de}
	(\e,\d)=
	(\frac{|\hx-\hy|}{4},\frac{\overline{C}^{-1}}{4}|\hx-\hy|^{\theta-2}).
\end{equation}
 With this pair of numbers, by taking $\tau>0$ small enough, there exists a constant $C>M$ such that the following inequalities hold.\\
(a) 
\begin{equation}\label{lem231}
	F(\hx,p+\tau \hx,Q+\tau)-F(\hy,p,-Q)
	\geq -C(|\hx-\hy|^{\rho_1}|p|^q+|\hx-\hy|^{\rho_2}||Q||).
\end{equation}
(b)
$$
	\int_{|z|\geq \e} [v(\hx+z)-v(\hx)-{\bf 1_{|z|\leq 1}}(p+\tau \hx)z]
 	c(z)dz -\int_{|z|\geq \e} [v(\hy+z)-v(\hy)
$$
\begin{equation}\label{lem233}
-{\bf 1_{|z|\leq 1}}pz]
 	c(z)dz 
 	\leq C \tau^{\frac{1}{2}}|\hx-\hy|^{-\gamma}.
\end{equation}
\end{lemma}
$Proof$ $of$ $Lemma$ $3.4.$ By putting  $\hat{z}=\hx-\hy$ in Lemma 3.3, for 
$\d=\frac{\overline{C}^{-1}}{4}|\hx-\hy|^{\theta-2}$, we can take 
$$
	\e=\min\{\d \overline{C}|\hx-\hy|^{3-\theta}, \frac{1}{2}|\hx-\hy|\}
	=\frac{1}{4}|\hx-\hy|,
$$
 so that (\ref{vdev1}), (\ref{vdev2}) hold. \\
 (a) 
From the continuity of $F$, (\ref{delliptic}), and (\ref{Fcont}), since $Q\leq O$, for $r_0>0$ ($|\hx-\hy|<r_0$) small enough, 
$$
	F(\hx,p+\tau \hx,Q+\tau)-F(\hy,p,-Q)
	\qquad\qquad\qquad\qquad\qquad\qquad\qquad\qquad\qquad
$$
$$
	=F(\hx,p,Q)-F(\hx,p,-Q)+F(\hx,p,-Q)-F(\hy,p,-Q)+o(\tau)
	\qquad\qquad\qquad
$$
\begin{equation}\label{Festimate}
	\geq 
	-w(\hx-\hy)|p|^q - \eta(\hx-\hy)||Q||+o(\tau)
	\geq -C(|\hx-\hy|^{\rho_1}|p|^q+|\hx-\hy|^{\rho_2}||Q||),
\end{equation}
where $C>M$ is a constant.\\
(b) Since 
$\Phi(\hx,\hy)=v(\hx)-v(\hy)-C_{\theta}|\hx-\hy|^{\theta}
	-\frac{\tau}{2}\hx^2
	\geq \Phi(0,0)=0,$ from (\ref{vbound}), 
\begin{equation}\label{rhx}
	\frac{\tau}{2}\hx^2\leq 2M. 
\end{equation}
Thus, for $\tau\in (0,1)$
$$
	v(\hx+z)-v(\hy+z)-(v(\hx)-v(\hy))\leq 
	\frac{\tau}{2}(\hx+z)^2-\frac{\tau}{2}\hx^2\leq 
	\tau^{\frac{1}{2}}(2M|z|+|z|^2),
$$
and from this 
$$
	\int_{|z|\geq \e}[v(\hx+z)-v(\hx)-{\bf 1}_{|z|\leq 1}(p+\tau \hx)z]
	c(z)dz
	-\int_{|z|\geq \e} [v(\hy+z)-v(\hy)
$$
\begin{equation}\label{this}
	-{\bf 1}_{|z|\leq 1}pz]
 	c(z)dz
 	\leq 
 	\int_{|z|\geq \e} \tau^{\frac{1}{2}}(3M|z|+z^2)c(z)dz
 	\qquad\qquad\qquad\qquad\qquad\qquad\qquad\qquad
\end{equation}

From (\ref{ccont}) and (\ref{23de}) 
\begin{equation}\label{r}
	\int_{|z|\geq \e} \tau^{\frac{1}{2}} (3M|z|+z^2) c(z)dz
	\leq \tau^{\frac{1}{2}}C\max\{1,|\hx-\hy|^{1-\gamma}\}
	\leq
	C\tau^{\frac{1}{2}}|\hx-\hy|^{-\gamma},
\end{equation}
where $C>M$ is the constant. 
By plugging (\ref{r}) into (\ref{this}), we get (\ref{lem233}).\\

	We put  the estimates (\ref{lem231})-(\ref{lem233}) in (\ref{subst1}), and since $\nu>0$ can be taken arbitrarily small,
$$
	C^{-1}{C_{\theta}}|\hx-\hy|^{\theta-\gamma} 
	\leq C( |\hx-\hy|^{\rho_1}|p|^q+|\hx-\hy|^{\rho_2}||Q||
	\qquad\qquad\qquad
$$
\begin{equation}\label{hen}
	+ 2\tau^{\frac{1}{2}}|\hx-\hy|^{-\gamma} + M|\hx-\hy|^{\theta_0}).
\end{equation}
From (\ref{ctheta}), $\theta\in (0,\min\{1,\theta_0+\gamma\})$, 
$\rho_1+\gamma>q$, $\rho_2+\gamma>2$, 
and since we can take $\tau \in (0,1)$ arbitrarily small, for $r_0>0$ ($|\hx-\hy|<r_0$) small enough, we get a contradiction. Thus, the claim in Theorem 3.1 is proved. \\

$Proof$ $of$ $Theorem$ $3.2.$ We use the similar contradiction argument as in the proof of Theorem 3.1. For an arbitrary fixed number $\theta\in (0,1)$, and for  $r_0>0$ small enough,  let  $C_{\theta}>0$ be such that (\ref{ctheta}): 
$$
	C_{\theta}r^{\theta}= 2M, 
$$
 and we shall prove (\ref{vhold}) by the contradiction's argument. 
As before, 
assume that there exist $x',y'\in {\bf R^N}$ ($|x'-y'|<r_0$) such that 
$$
	v(x')-v(y')>C_{\theta}|x'-y'|^{\theta},
$$
and we shall look for a contradiction. However, we must modify the preceding argument,  because for $N\geq 2$ the matrix 
$Q=\n^2_{xx} \phi(\hx,\hy)=\n^2_{yy} \phi(\hx,\hy)$ ($\phi$ is the function in the proof of Theorem 3.1) is no longer negatively 
definite, and we have to use  Lemma 2.4. For this reason, let us consider the supconvolution $v^r$ and the infconvolution $v_r$ of $v$  defined by (\ref{supconvolution}) and (\ref{infconvolution}), respectively. 
From Lemma 2.3, for any $\nu>0$ there exists $r_1>0$ such that $v^r$ ($\forall r\in (0,r_1)$) is a subsolution of 
$$
	F(x,\n v^r(x),\n^2 v^r(x))-\int_{{\bf R^N}}  
	[v^r(x+z)-v^r(x)
		\qquad\qquad\qquad
$$
\begin{equation}\label{supfirst}
	-{\bf 1}_{|z|\leq 1}\la \n v^r(x),z \ra]c(z) dz -g(x) \leq \nu \qquad x\in {\bf R^N},
\end{equation}
and $v_r$ ($\forall r\in (0,r_1)$) is a supersolution of 
$$
	F(x,\n v_r(x),\n^2 v_r(x))-\int_{{\bf R^N}}  
	[v_r(x+z)-v_r(x)
		\qquad\qquad\qquad
$$
\begin{equation}\label{supfirst}
	-{\bf 1}_{|z|\leq 1}\la \n v_r(x),z \ra]c(z) dz -g(x) \geq -\nu \qquad x\in {\bf R^N}. 
\end{equation}
Of course from the preceding assumption, for $\forall r\in (0,r_0)$
$$
	v^r(x')-v_r(y')>C_{\theta}|x'-y'|^{\theta}. 
$$
Now, consider for $\tau\in (0,1)$ 
$$
	\Phi(x,y)=v^r(x)-v_r(y)-C_{\theta}|x-y|^{\theta}-\frac{\tau}{2}|x|^2,
$$
and let $(\hat{x},\hat{y})$ be a maximum point of $\Phi$.
 Put $p=\n_x \phi(\hx,\hy)=-\n_y \phi(\hx,\hy)$, and $Q=\n^2_{xx} \phi(\hx,\hy)=\n^2_{yy} \phi(\hx,\hy)$. We  use Lemma 2.4 for $U=v^{r}-\frac{\tau}{2}|x|^2$, and $V=v_r$, and for $O=\{(x,y)\in {\bf R^{2N}}| |x-y| < r_0 \}$, and 
 we know that there exists $(x_m,y_m)\in {\bf R^{2N}}$ such that 
 $\lim_{m\to \infty}(x_m,y_m)=(\hat{x},\hat{y})$. There also exist 
 $(p_m+\tau x_m,X_m+\tau I)\in J^{2,+}_{{\bf R^N}}v^r(x_m)$,  $(p_m',Y_m)\in J^{2,-}_{{\bf R^N}}v_r(y_m)$
 such that 
$\lim_{m\to \infty}p_m$$=\lim_{m\to \infty}p'_m$$=2\a(x_m-y_m)=p$, and 
$X_m\leq Y_m\quad \forall m$. Moreover, the claim in Lemma 2.4 (iii) leads 
 the following  for any $z\in {\bf R^N}$ such that $(x_m+z,y_m+z)\in O$
\begin{equation}\label{important2}
	v^r(x_m+z)-v^r(x_m)-\la p_m,z\ra
	-\{ v_r(y_m+z)-v_r(y_m)-\la p'_m,z\ra\} 
\end{equation}
$$
	\leq \frac{\tau}{2}|x_m+z|^2-\frac{\tau}{2}|x_m|^2
	=\frac{\tau}{2}\{2\la x_m,z\ra+|z|^2\}.
$$
Let 
$(\e_m,\d_m)$ be a pair of positive numbers such that 
\begin{equation}\label{mm}
	v^r(x_m+z)\leq v^r({x}_m)+\la (p_m+\tau x_m),z \ra 
	+\frac{1}{2}\la (X_m+\tau I) z,z \ra 
	+\delta_m |z|^2  \qquad \hbox{if}\quad |z|\leq \e_m,
\end{equation}
and
\begin{equation}\label{MM}
	v_r(y_m+z)\geq v_r({y}_m)+\la p'_m,z \ra 
	+\frac{1}{2}\la Y_m z,z \ra 
	-\delta_m |z|^2  \qquad \hbox{if}\quad |z|\leq \e_m.\qquad\qquad\qquad
\end{equation}
Then, from the definition of viscosity solutions, we have 
$$
	F(x_m,p_m+\tau x_m,X_m+\tau I)
 	-\int_{|z|< \e_m}
 	\frac{1}{2}\la (X_m+(\tau+2\d_m)I)z,z\ra c(z)dz 	 
$$
$$
 	-\int_{|z|\geq \e_m} [v^r(x_m+z)-v^r(x_m)
 	-{\bf 1}_{|z|\leq 1}\la (p_m+\tau x_m),z\ra]
 	c(z)dz-g(x_m) \leq \nu,
$$
 and 
 $$
 	F(y_m,p_m,Y_m)
 	-\int_{|z|< \e_m}\frac{1}{2}\la (Y_m-2\d_m I)z,z \ra c(z)dz
$$
$$
	-\int_{|z|\geq \e_m} [v_r(y_m+z)-v_r(y_m)
	-{\bf 1}_{|z|\leq 1}\la p_m',z\ra ]
 	c(z)dz 
 	-g(y_m) \geq -\nu.
$$
By taking the difference of the two inequalities, 
$$
	F(x_m,p_m+\tau x_m,X_m+\tau I)-F(y_m,p_m,Y_m)
	\qquad\qquad\qquad\qquad\qquad\qquad\qquad\qquad\qquad
$$
$$
	-\frac{1}{2}\int_{|z|\leq \e_m}\la (X_m-Y_m+(\tau+4\d_m)I)z,z \ra 
	c(z)dz
$$
$$
	\leq 2 \nu + 
	\int_{|z|\geq \e_m} [v^r(x_m+z)-v^r(x_m)
	-{\bf 1}_{|z|\leq 1}\la p_m+\tau x_m, z \ra]
 	c(z)dz 
 	\qquad\qquad\qquad\qquad\qquad\qquad
$$
$$
		-\int_{|z|\geq \e_m} [v_r(y_m+z)-v_r(y_m)
		-{\bf 1}_{|z|\leq 1}\la p'_m,z \ra]
 	c(z)dz +g(x_m)-g(y_m)
$$
$$
	\leq 
	2\nu + \int_{{|z|\geq \e_m}\cap O_m(z)}\frac{\tau}{2}|z|^2 c(z)dz
	\int_{{|z|\geq \e_m}\cap O_m(z)^c} 
	[\{v^r(x_m+z)-v^r(x_m)
$$
$$
	-{\bf 1}_{|z|\leq 1}\la p_m+\tau x_m, z \ra\}
	-\{v_r(y_m+z)-v_r(y_m)
		-{\bf 1}_{|z|\leq 1}\la p'_m,z \ra\}
	],
$$
where
$$
	O_m(z)=\{
	z\in {\bf R^N}|\quad (x_m+z,y_m+z)\in \Omega\times \Omega
	\}.
$$
Since $\lim_{m\to \infty}(x_m,y_m)=(\hat{x},\hat{y})\in \Omega\times \Omega$,  there exists a ball $B(0)\subset {\bf R^N}$, centered at the origin, independent on $m\in {\bf N}$, such that $B(0)\subset O(z)=\lim_{m\to \infty}O_m(z)$, i.e.
\begin{equation}\label{ball}
	(x_m+z,y_m+z)\in \Omega\times \Omega\quad \forall z\in B(0),\quad \forall m\in {\bf N}.
\end{equation}
Then, by pussing $m\to \infty$ in the above inequality, we get
$$
	F(\hat{x},p+\tau \hat{x},X+\tau I)-F(\hat{y},p,Y)
	\qquad\qquad\qquad\qquad\qquad\qquad\qquad\qquad\qquad
$$
$$
	\leq 
	2\nu + \int_{O(z)}\frac{\tau}{2}|z|^2 c(z)dz + 
	\int_{O(z)^c} 
	[\{v^r(\hat{x}+z)-v^r(\hat{x})
$$
$$
	-{\bf 1}_{|z|\leq 1}\la p+\tau \hat{x}, z \ra\}
	-\{v_r(\hat{y}+z)-v_r(\hat{y})
		-{\bf 1}_{|z|\leq 1}\la p,z \ra\}
	]c(z)dz.
$$
$$
	\leq C(\nu +M)+ \int_{{\bf R^N}}\frac{\tau}{2}|z|^2 c(z)dz
	+ \int_{\{|z|<1\}\cap O(z)^c} \tau|\hat{x}||z| c(z)dz,
$$
where $C>0$ is a constant, and we have used the fact that $(\hat{x},\hat{y})$ 
is the maximizer of $\Phi$. From (\ref{cintegral}) and (\ref{ctheta}), and since $O(z)^c \subset B(0)^c$, for $0<\tau<1$, 
\begin{equation}\label{subst2}	
	F(\hat{x},p+\tau \hat{x},X+\tau I)-F(\hat{y},p,Y)
	\leq C(\nu +M+\tau^{\frac{1}{2}}),
\end{equation}
where $C>0$ is a constant. We shall give the estimate of the left-hand side of the above.\\

{\bf Lemma 3.5.}
\begin{lemma}
There exists a constant $C>M$ such that the following holds.\\
\begin{equation}\label{lem241}
	F(\hx,p+\tau \hx,X+\tau I)-F(\hy,p,-Y)\geq \frac{C_{\theta}}{C}|\hx-\hy|^{\theta-2} + o(\tau).
\end{equation}
\end{lemma}

$Proof$ $of$ $Lemma$ $3.5.$ 
 From the continuity of $F$, (\ref{Fcont}), (\ref{uelliptic}), and (\ref{rhx}) (which is also true for $N\geq 2$), 
$$
	F(\hx,p+\tau \hx,X+\tau I)-F(\hy,p,-Y)= 
	F(\hx,p,\overline{X})-F(\hy,p,-\overline{Y})+o(\tau)+o(\nu')
$$
$$
	= F(\hx,p,\overline{X})-F(\hx,p,-\overline{Y})
	+	F(\hx,p,-\overline{Y})-F(\hy,p,-\overline{Y})+o(\tau)+o(\nu')
$$
\begin{equation}\label{Funif}
	\geq -\l_0 \hbox{Tr}(\overline{X}+\overline{Y})-w(|\hx-\hy|)|p|^q-\eta(|\hx-\hy|)||\overline{Y}||+o(\tau)+o(\nu').
\end{equation}

We need the following lemma, the proof of which is delayed in the end.\\
{\bf Lemma 3.6.}
\begin{theorem} If $A$, $B$, and $Q\in {\bf S^N}$ satisfy 
\begin{equation}\label{mineq}
	\left(
		\begin{array}{cc}
			A&O\\
			O&B
		\end{array}
	\right) \leq 
		\left(
		\begin{array}{cc}
			Q &-Q\\
			-Q&Q
		\end{array}
	\right),
\end{equation}
then there exists a constant $L>0$ such that 
$$
	||A||,\quad ||B||\quad \leq 
	L||Q||^{\frac{1}{2}}|\hbox{Tr}(A+B)|^{\frac{1}{2}}.
$$
The constant $L$ depends only on $N$.
\end{theorem}

Remark that 
\begin{equation}\label{XY}
	-Tr(X-Y)\geq C_{\theta}\theta(1-\theta)|\hx-\hy|^{\theta-2}>0,
\end{equation}
 because , 
$X-Y\leq 2Q$, $X-Y\leq O$, and for 
$O \leq P=\frac{(\hx-\hy)\otimes(\hx-\hy)}{|\hx-\hy|^2}\leq I$, 
$$
	\hbox{Tr}({X}-{Y})\leq 
	\hbox{Tr}(P({X}-{Y}))
	\leq 2\hbox{Tr}(PQ)=2C_{\theta}\theta(\theta-1)|\hx-\hy|^{\theta-2}<0.
$$
Therefore, by putting $A={X}$ and $B={Y}$ in Lemma 3.6, and by taking $r_0>0$ ($|\hx-\hy|<r_0$) small enough, 
from (\ref{XY}) 
$$
	\eta(|\hx-\hy|)||{Y}||
	\leq K'C_{\theta}\eta(|\hx-\hy|)|\hx-\hy|^{\theta-2}
	\leq KC_{\theta}|\hx-\hy|^{\theta-2},
$$
where $K,K'>0$ are constants. 
For $r_0>0$ ($|\hx-\hy|<r_0$) small enough, from (\ref{os211}) and (\ref{ctheta})
$$
	w(|\hx-\hy|)|p|^q= w(|\hx-\hy|)C_{\theta}^q |\hx-\hy|^{q(\theta-1)}
	\qquad\qquad\qquad\quad
$$
$$
	\qquad\qquad\qquad
	\leq L(C_{\theta}|\hx-\hy|^{\theta})^q|\hx-\hy|^{\rho_1-q}
	\leq \frac{\l_0}{4}C_{\theta}|\hx-\hy|^{\theta-2}. 
$$
Therefore, from (\ref{Funif}) and (\ref{XY}), 
$$
	F(\hx,p+\tau \hx,X+\tau I)-F(\hy,p,Y)\geq
	 \frac{C_{\theta}}{C}|\hx-\hy|^{\theta-2} +o(\tau),\quad
$$
where $C>M$ is a constant. We showed (\ref{lem241}).\\

By plugging (\ref{lem241}) into (\ref{subst2}), since $\nu>0$ can be taken arbitrarily small, for any $0<\theta < 1$, 
we get a contradiction for $r_0>0$ ($|\hx-\hy|<r_0$) small enough. 
We have proved (\ref{vhold}). \\

Finally, we are to prove Lemma 3.6. \\

$Proof$ $of$ $Lemma$ $3.6$ By multiplying the matrix 
$$
	\left(
		\begin{array}{cc}
			I&I\\
			I&-I
		\end{array}
	\right) 
$$
to the both hand sides of (\ref{mineq}) first from right and then from left, we get
$$
	\left(
		\begin{array}{cc}
			A+B&A-B\\
			A-B&A+B
		\end{array}
	\right) 
	\leq
	\left(
		\begin{array}{cc}
			O&O\\
			O&4Q
		\end{array}
	\right).
$$
	Thus, for any $t\in {\bf R}$ and $\xi\in {\bf R^N}$
$$
	\left(
		\begin{array}{cc}
		t\xi&\xi\\
		\end{array}
	\right)
	\left(
		\begin{array}{cc}
			A+B&A-B\\
			A-B&A+B
		\end{array}
	\right) 
	\left(
		\begin{array}{c}
			t\xi\\
			\xi
		\end{array}
	\right)
	\leq
	\left(
		\begin{array}{cc}
		t\xi&\xi\\
		\end{array}
	\right)
	\left(
		\begin{array}{cc}
			O&O\\
			O&4Q
		\end{array}
	\right)
	\left(
		\begin{array}{c}
			t\xi\\
			\xi
		\end{array}
	\right),
$$
and 
$$
	t^2\la \xi,(A+B)\xi\ra +2t\la \xi,(A-B)\xi\ra +\la \xi,(A+B)\xi\ra 
	-4\la \xi,Q\xi\ra \leq 0.
$$
Hence, for any $|\xi|=1$, 
$$
	\la \xi,(A-B)\xi\ra^2 \leq 
	\la \xi,(A+B)\xi\ra (4\la \xi,Q\xi\ra-\la \xi,(A+B)\xi\ra).
$$
This yields $\la \xi,(A+B)\xi\ra^2 \leq 4||A+B||\cdot||Q||$, and since 
$||A+B||\leq C|\hbox{Tr}(A+B)|\cdot||Q||$ where $C>0$ is a constant depending only on  $N>0$, we proved the claim.\\

%%%%%%%%%%%%%%%%%%%%%%%%%%%%%%%%%%%%%%%%%%%%%%%%%%%%%%%%%%%%%%%%%%%%%%%%%%%%%%
\section{Other H{\"o}lder continuities of viscosity solutions}

	In this section, we shall study (\ref{lambH}): 
$$
	\l v(x)+ H(\n v(x))-\int_{{\bf R^N}}  
	[v(x+z)-v(x)
	\qquad\qquad\qquad\qquad\qquad\qquad
$$
$$
	-{\bf 1}_{|z|\leq 1}\la \n v(x),z \ra]c(z) dz -g(x)=0 \qquad x\in 
	{\bf T^N},
$$
and  (\ref{lamb}): 
$$
	\l v(x)+ F(x,\n v(x),\n^2 v(x)) -\int_{{\bf R^N}}  
	[v(x+z)-v(x)
	\qquad\qquad\qquad\qquad
$$
$$
	-{\bf 1}_{|z|\leq 1}\la \n v(x),z \ra]c(z) dz -g(x)=0 \qquad x\in {\bf T^N},
$$
where $\l>0$. We consider the case other than 
(I) $N=1$, and (II) $F$ is uniformly elliptic. So, we are interested in  the case of $N\geq 2$, and $F$ (or $H$) is degenerate elliptic. We assume the conditions (\ref{periodicity})-(\ref{pFcont}).
\\

{\bf Example 4.1.\quad} The following is an 
example of $F$ satisfying the conditions (\ref{periodicity})-(\ref{pFcont}).$$
	-\sum_{i=1}^{N-1} a_i(x)\frac{\p^2 u}{\p x_i^2}(x) 
	-\sum_{i=1}^{N-1} b_i(x)\frac{\p u}{\p x_i}(x)
	+|\frac{\p u}{\p x_N}(x)| \qquad x\in {\bf T^N},
$$
where $a_i(x)>\exists \l_1>0$ and $b_i(x)$ ($1\leq i\leq N-1$) are periodic in ${\bf T^N}$. 
Or, more generally the following Hamilton-Jacobi-Bellman operator satisfies 
(\ref{periodicity})-(\ref{pFcont}).
$$
	F(x,u,\nabla u,\nabla^2 u)=\sup_{\alpha \in \CA} \{ 
	-\sum_{ij=1}^N a_{ij}(x,\alpha)\frac{\p^2 u}{\p x_i \p x_j}
	-\sum_{i=1}^N b_{i}(x,\alpha)\frac{\p u}{\p x_i}\qquad \qquad \qquad \qquad
$$
\begin{equation}\label{a3}	
	\qquad \qquad \qquad\qquad \qquad \qquad \qquad 
	+c(x,\alpha)u-f(x,\alpha)\} \qquad x\in {\bf T^N},
\end{equation}	
where $\CA$ a given set (controls), $(a_{ij}(x,\alpha)\in {\bf S}^N$ ($\alpha \in \CA$) non-negative matrices periodic in ${\bf T^N}$ such that there exist  matrices 
 $\sigma^{\alpha}$ ($\alpha \in A$) of the size $N\times k$, 
$$
	A_{\a}=(a_{ij}(x,\alpha))=\sigma(x,\a)^T\sigma(x,\a),	
$$
$$
	A_{\a}'\geq \l_1 I_M,
	\quad A_{\a}=\left(
		\begin{array}{cc}
			A_{\a}'&{A_{\a}}_{12}\\
			{A_{\a}}_{21}&{A_{\a}}_{22}
		\end{array}
	\right), \quad A_{\a}'\in {\bf S}^M \quad (M<N),
$$
where $\l_1>0$, 
and 
$b(x,\alpha)=$$(b_i(x,\alpha))\in {\bf R}^N$, $c(x,\a)\in {\bf R}$
 are bounded, periodic in  ${\bf T^N}$, and regular enough. \\

 We shall give the results.\\
 
{\bf Theorem 4.1.\quad}  
\begin{theorem} Let $v$ be a periodic viscosity solution of (\ref{lambH}) satisfying (\ref{vbound}). Assume that (\ref{cintegral}), (\ref{ccont}), (\ref{gcont}), and (\ref{periodicity}) hold, where $\gamma\in (0,2)$, and $\l>0$. Let  $H(p)=|p|^q$, where  $q\geq 1$.
 Then, for any $\theta\in (0,\theta_0)$, there exists a constant $C_{\theta}>0$ such that (\ref{lambvhold}) holds. 
 The constant $C_{\theta}$ does not depend on $\lambda\in (0,1)$.
 \end{theorem}
 
{\bf Theorem 4.2.\quad}  
\begin{theorem}  Let $v$ be a periodic viscosity solution of (\ref{lamb}) satisfying (\ref{vbound}). Assume that (\ref{cintegral}), (\ref{ccont}), (\ref{gcont}), (\ref{periodicity}), (\ref{A1}) and (\ref{pFcont}) hold, where $\gamma\in (0,2)$, $\l>0$, and that 
 there exist constants $L>0$, $\rho_i>0$ ($i=1,2$) such that 
\begin{equation}\label{os321}
	\lim_{s\downarrow 0}w(s)s^{-\rho_1}\leq L,\quad 
	\lim_{s\downarrow 0}\eta(s)s^{-\rho_2}\leq L,
\end{equation}
where $\rho_1+\gamma>q$, $\rho_2+\gamma>2$. 
 Then, for any $\theta\in (0,\theta_0)$ there exists a constant $C_{\theta}>0$ such that (\ref{lambvhold}) holds. 
 The constant $C_{\theta}$ does not depend on $\lambda\in (0,1)$.
\end{theorem}

$Proof$ $of$ $Theorem$ $4.1.$ 
We use the contradiction argument similar to that of Theorem 2.1. Fix $\theta\in (0,\theta_0)$, and let $r_0>0$ be small enough. Let us  take 
$\overline{C}_{\theta}>0$ such that 
\begin{equation}\label{lctheta}
	\overline{C}_{\theta}r^{\theta}= 2M, 
\end{equation}
and we shall prove (\ref{lambvhold}) (for $\frac{C_{\theta}}{\l}=\overline{C_{\theta}}$) by contradiction. 
For $x,y\in {\bf T^N}$ such that $|x-y|\geq r_0$, from (\ref{vbound}) we have 
$$
	|v(x)-v(y)|\leq 2M\leq \overline{C}_{\theta}|x-y|^{\theta}.
$$
Thus, assume that there exist $x',y'\in {\bf T^N}$ ($|x'-y'|<r_0$) such that 
$$
	v(x')-v(y')> \overline{C}_{\theta}|x'-y'|^{\theta},
$$
and we shall look for a contradiction. As in the proof of Theorem 3.2,  take the supconvolution $v^r$ and the infconvolution $v_r$ of $v$, which are respectively the subsolution and the supersolution of the following problems. 
$$
	\l v^r(x)+ H(\n v^r(x))-\int_{{\bf R^N}}  
	[v^r(x+z)-v^r(x)
	\qquad\qquad\qquad\qquad\qquad\qquad
$$
$$
	-{\bf 1}_{|z|\leq 1}\la \n v^r(x),z \ra]c(z) dz -g(x)\leq \nu \qquad x\in 
	{\bf T^N},
$$
$$
	\l v_r(x)+ H(\n v_r(x))-\int_{{\bf R^N}}  
	[v_r(x+z)-v_r(x)
	\qquad\qquad\qquad\qquad\qquad\qquad
$$
$$
	-{\bf 1}_{|z|\leq 1}\la \n v_r(x),z \ra]c(z) dz -g(x)\geq -\nu \qquad x\in 
	{\bf T^N},
$$
where $\nu>0$ is an arbitrary small constant. Remark that 
$$
	v^r(x')-v_r(y')> \overline{C}_{\theta}|x'-y'|^{\theta}
$$
holds for $|x'-y'|<r_0$. 
Consider 
$$
	\Phi(x,y)=v^r(x)-v_r(y)-\overline{C}_{\theta}|x-y|^{\theta},
$$
and let $(\hat{x},\hat{y})$ be a maximum point of $\Phi$. Let us write $\phi(x,y)=\overline{C}_{\theta}|x-y|^{\theta}$. For
$$
	\n_x \phi(x,y)=\overline{C}_{\theta} \theta |x-y|^{\theta-2}(x-y)=-\n_y \phi(x,y), 
$$
$$
	\n^2_{xx} \phi(x,y)=\overline{C}_{\theta}\theta |x-y|^{\theta-2}I+
	\overline{C}_{\theta}\theta(\theta-2) |x-y|^{\theta-2}(x-y)\otimes(x-y)
	=\n^2_{yy} \phi(x,y), 
$$
put $p=\n_x \phi(\hx,\hy)=-\n_y \phi(\hx,\hy)$, and $Q=\n^2_{xx} \phi(\hx,\hy)=\n^2_{yy} \phi(\hx,\hy)$. As in the proof of Theorem 3.2, by using Lemma 2.4 for $U=v^{r}$, $V=v_r$, and $O=\{(x,y)\in {\bf R^{2N}}| |x-y| < r_0 \}$, 
 we know that there exists $(x_m,y_m)\in {\bf T^{2N}}$ such that 
 $\lim_{m\to \infty}(x_m,y_m)=(\hat{x},\hat{y})$. There also exist 
 $(p_m,X_m)\in J^{2,+}_{{\bf T^N}}v^r(x_m)$,  $(p_m',Y_m)\in J^{2,-}_{{\bf T^N}}v_r(y_m)$
 such that 
$\lim_{m\to \infty}p_m$$=\lim_{m\to \infty}p'_m$$=2\a(x_m-y_m)=p$, and 
$X_m\leq Y_m\quad \forall m$. The claim in Lemma 2.4 (iii) leads 
  for any $z\in {\bf R^N}$ such that $(x_m+z,y_m+z)\in O$, 
\begin{equation}\label{important3}
	v^r(x_m+z)-v^r(x_m)-\la p_m,z\ra
	-\{ v_r(y_m+z)-v_r(y_m)-\la p'_m,z\ra\} \leq 0. 
\end{equation}
Let 
$(\e_m,\d_m)$ be a pair of positive numbers such that 
\begin{equation}\label{mm3}
	v^r(x_m+z)\leq v^r({x}_m)+\la p_m,z \ra 
	+\frac{1}{2}\la X_m z,z \ra 
	+\delta_m |z|^2  \qquad \hbox{if}\quad |z|\leq \e_m,
\end{equation}
and
\begin{equation}\label{MM3}
	v_r(y_m+z)\geq v_r({y}_m)+\la p'_m,z \ra 
	+\frac{1}{2}\la Y_m z,z \ra 
	-\delta_m |z|^2  \qquad \hbox{if}\quad |z|\leq \e_m.
\end{equation}

By using the similar argument as in Theorem 3.2, from the definition of viscosity solutions, we have the following. 
$$
	\l (v^r(x_m)-v_r(y_m))+ H(p_m)-H(p'_m)\qquad\qquad\qquad\qquad
	\qquad\qquad\qquad\qquad\qquad
$$
$$
	-\frac{1}{2}\int_{|z|\leq \e_m}\la (X_m+2\d_m I)z,z \ra 
	c(z)dz
	-\frac{1}{2}\int_{|z|\leq \e_m}\la (Y_m-2\d_m I)z,z \ra  
	c(z)dz\qquad
$$
$$
	-\int_{|z|\geq \e_m} [v^r(x_m+z)-v^r(x_m)
	-{\bf 1}_{|z|\leq 1}\la p_m,z \ra]
 	c(z)dz 
 	\qquad\qquad\qquad\qquad\qquad\qquad
$$
$$
	+\int_{|z|\geq \e_m} [v_r(y_m+z)-v_r(y_m)
	-{\bf 1}_{|z|\leq 1}\la p'_m,z \ra]
 	c(z)dz 
 	\leq g(x_m)-g(y_m)+2\nu.
$$
Remarking that $X_m\leq Y_m$ and (\ref{important3}) hold, and 
by using the similar argument as in Theorem 3.2, we can pass $m\to \infty$ in 
the above inequality to have 
$$
	\l (v^r(\hat{x})-v_r(\hat{y})) \leq 
	g(\hat{x})-g(\hat{y})+2\nu,
$$
and since $\nu>0$ is arbitrary, 
we have 
$$
	\l \overline{C}_{\theta}|\hx-\hy|^{\theta}\leq M|\hx-\hy|^{\theta_0}.
$$
Since $\overline{C}_{\theta}=\frac{2M}{r_0^{\theta}}$, the above leads 
$$
	2\l\leq r_0^{\theta_0}.  
$$
However, if we take for an arbitrarily fixed $c>\frac{1}{\theta_0}$, 
\begin{equation}\label{C}
	r_0= \l^{c},\quad \overline{C}_{\theta}=\frac{2M}{\l^{c\theta}},
\end{equation}
we get a contradiction for any $\l\in (0,1)$. Therefore, for $0<\theta<\theta_0$, by taking $c=\frac{1}{\theta}$  and thus $\overline{C}_{\theta}=\frac{2M}{\l}$, we proved our claim for $C_{\theta}=2M$
$$
	v(x)-v(y)\leq \overline{C}_{\theta}|x-y|^{\theta}=
	\frac{2M}{\l}|\hx-\hy|^{\theta} \quad \forall x,y\in {\bf T^N}.
$$
$$\quad$$

$Proof$ $of$ $Theorem$ $4.2.$ The argument is similar to that of Theorem 4.1, 
 and we omit the proof. \\

%%%%%%%%%%%%%%%%%%%%%%%%%%%%%%%%%%%%%%%%%%%%%%%%%%%%%%%%%%%%%%%%%%%%%%%%%%%%%
\section{Strong maximum principle}

$\quad$In this section, we consider 
\begin{equation}\label{strong}
	F(x,\n u,\n^2 u) -\int_{{\bf R^N}}  
	[u(x+z)-u(x)
	-\la \n u(x),z \ra] c(z) =0 
	\qquad \forall x\in {\bf R^N},
\end{equation}
where $F$ satisfies (\ref{delliptic}) and 
\begin{equation}\label{positive}
	F(x,0,O)\quad \geq \quad 0\quad \forall x\in {\bf R^N}.
\end{equation}
	We  assume the following condition.
	
\leftline{(Almost everywhere positivity) : For 
any open set $D\in {\bf R^N}$,}
\begin{equation}\label{aepositive}
	\int_{z\in D} 1 c(z)dz >0. 
\end{equation}
$$\quad$$
Our strong maximum principle is the following.\\
{\bf Theorem 5.1 (\cite{mpl2}).}
\begin{theorem} Consider the integro-differential equation (\ref{strong}), and 
 assume that (\ref{cintegral}), (\ref{ccont}), (\ref{positive}), and (\ref{aepositive}) hold. 
  Let $u$ be a viscosity subsolution of 
  (\ref{strong}), and assume that it takes a maximum at a point $x_0\in {\bf R^N}$, i.e. 
\begin{equation}\label{glob}
	u(x) \leq u(x_0) \qquad \forall x\in {\bf R^N}.
\end{equation}
 Then, $u$ is constant in ${\bf R^N}$ almost everywhere.
\end{theorem}

$Proof.$ $\quad$ From (\ref{glob}), for $p={0}$ and $X=O$, 
$$
	u(x_0+z)\leq u(x_0)+\la 0, z\ra+\frac{1}{2}\la Oz,z\ra +\delta |z|^2
	\quad \hbox{if}\quad |z|\leq \e
$$
holds for any $\d>0$ and $\e>0$. Hence, from the definition of viscosity 
subsolution
$$
	F(x_0,0,O)-\int_{|z|\leq \e}  
	\frac{1}{2} \la (O+2\d I)z,z \ra c(z)dz
$$
$$
	-\int_{|z|\geq \e} [u(x_0+z)
	-u(x_0)-\la 0,z \ra] c(z)\}\leq 0, 
$$
holds for any $\d>0$ and $\e>0$. So, from (\ref{positive}) we have 
$$
	\int_{|z|\geq \e} [u(x_0)-u(x_0+z)] c(z) \leq 0
$$
holds for any $\e>0$. Therefore,  from (\ref{ccont}), 
  (\ref{aepositive}), and (\ref{glob}), 
$$
	u(x)\leq u(x_0)\leq u(x) \qquad \hbox{almost everywhere in}\quad 
	x\in {\bf R^N},
$$
and the claim is proved. \\

{\bf Remark 5.1.} We shall use the above strong maximum principle to solve the 
ergodic problem in the next section. 
$$\quad$$

%%%%%%%%%%%%%%%%%%%%%%%%%%%%%%%%%%%%%%%%%%%%%%%%%%%%%%%%%%%%%%%%%%%%%%%%%%%
\section{Ergodic problem for integro-differential equations}

$\quad$In this section, we apply the results in preceding sections to solve the ergodic 
 problem in ${\bf T^N}$. We shall study the existence of a unique number $d_f$ such that the following problem has a periodic viscosity solution.
$$
	d_f+F(x,\n u,\n^2 u)-\int_{{\bf R^N}}  
	[u(x+z)-u(x)
	\qquad\qquad\qquad\qquad\qquad
$$
\begin{equation}\label{ergodic}
	-{\bf 1}_{|z|\leq 1}\la \n u(x),z \ra] c(z)dz -f(x)=0 \qquad x\in {\bf T^N}.
\end{equation}
 For this purpose, we consider  the approximated problem: 
\begin{equation}\label{ax}
	\l u_{\l}+F(x,\n u_{\l},\n^2 u_{\l})-\int_{{\bf R^N}}  
	[u_{\l}(x+z)-u_{\l}(x)
	\qquad\qquad\qquad\qquad\qquad
\end{equation}
$$
	-{\bf 1}_{|z|\leq 1}
	\la \n u_{\l}(x),z \ra] c(z)dz -f(x)=0 \qquad x\in {\bf T^N},
$$
where $\l\in (0,1)$, and we shall see whether there exists the following unique limit number 
$$
	\lim_{\l\downarrow 0} \l u_{\l}(x)=d_f \qquad \hbox{uniformly in}\quad {\bf T^N}.
$$
 We assume that $F$ satisfies (\ref{a2}), and that the following hold.\\
 
(Periodicity) :
\begin{equation}\label{periodic}
	F(\cdot,p,X),\quad f(\cdot) \quad \hbox{are periodic in}
	\quad x\in {\bf T^N},\quad \forall 
	(p,X)\in ({\bf R^N}\times {\bf S^N}). 
\end{equation}

(Homogeneity) : The partial differential operator $F$ is positively homogenious in degree one 
$$
	F(x,\xi p,\xi X)=\xi F(x,p,X)\qquad\qquad\qquad\qquad\qquad\qquad\qquad
$$
\begin{equation}\label{homog}
	\qquad\qquad\qquad\qquad
	 \quad \forall \xi >0,\quad \forall x\in {\bf T^N}, 
	 \quad \forall p\in {\bf R^N},
	 \quad \forall X\in {\bf S^N}.
\end{equation}
	
	As we have seen in  Theorem 2.2, under (\ref{a2}) and (\ref{periodic}), the comparison result holds.  From the Perron's method (see \cite{ar4} and \cite{ar5}), it is known that there exists a unique periodic viscosity solution $u_{\l}$ of (\ref{ax}) for any $\l\in (0,1)$. 
 Now, we state our main result.\\
	
{\bf Theorem 6.1.}
\begin{theorem}
	Let $u_{\l}$ ($\l\in (0,1)$) be the periodic viscosity solution of (\ref{ax}). Assume that the conditions in Theorem 5.1, (\ref{a2}), (\ref{periodic}), and (\ref{homog}) hold. Fix an arbitrary point $x_0\in {\bf T^N}$. Then, the following hold.\\
(i) Assume that the conditions in Theorem 3.1, or those in Theorem 3.2 hold. 
Then, there exist a unique number 
$d_f$ and a periodic function $u$ such that 
\begin{equation}\label{lconv}
	\lim_{\l\downarrow 0}\l u_{\l}(x)=d_f,\quad
	 \lim_{\l\downarrow 0} (u_{\l}(x)-u_{\l}(x_0))=u(x)
	 \quad\hbox{uniformly in}\quad {\bf T^N},
\end{equation}
such that (\ref{ergodic}) holds in the sense of viscosity solutions.\\ 
(ii) Assume that the conditions in Theorem 4.1, or those in Theorem 4.2 hold. Then, there exists a unique number 
$d_f$ such that
$$
	\lim_{\l\downarrow 0}\l u_{\l}(x)=d_f 
	\quad\hbox{uniformly in}\quad {\bf T^N},
$$
which is characterized by the following. For any $\nu>0$ there exist
a periodic viscosity subsolution $\underline{u}$ and a periodic viscosity supersolution $\overline{u}$ of 
$$
	d_f+F(x,\n \underline{u},\n^2 \underline{u})-\int_{{\bf R^N}}  
	[\underline{u}(x+z)
	\qquad\qquad\qquad\qquad\qquad\qquad\qquad
$$
$$
	\qquad\qquad
	-\underline{u}(x)-{\bf 1}_{|z|\leq 1}\la \n \underline{u}(x),z \ra] c(z)dz 
	-f(x)\leq \nu \qquad x\in {\bf T^N},
$$
$$
	d_f+F(x,\n \overline{u},\n^2 \overline{u}) -\int_{{\bf R^N}}  
	[\overline{u}(x+z)
	\qquad\qquad\qquad\qquad\qquad\qquad\qquad
$$
$$
	\qquad\qquad
	-\overline{u}(x)-{\bf 1}_{|z|\leq 1}\la \n \overline{u}(x),z \ra] c(z)dz -f(x)\geq -\nu \qquad x\in {\bf T^N}.
$$
\end{theorem}

$Proof$ $of$ $Theorem$ $6.1.$ (i) 
We shall prove the claim in the following three steps.\\
(Step 1.) We prove the uniform boundedness of $v_{\l}(x)=u_{\l}(x)-u_{\l}(x_0)$ :   
\begin{equation}\label{M}
	|v_{\l}(x)|=|u_{\l}(x)-u_{\l}(x_0)|\leq \exists M' \quad \forall x\in {\bf T^N},
	\quad \forall \l\in (0,1), 
\end{equation}
by a contradiction argument. Assume that there 
exists a subsequence $\l'\to 0 $ such that 
\begin{equation}\label{increase}
	\lim_{\l'\to 0}|v_{\l'}|_{L^{\infty}}=\infty,
\end{equation}
and we shall look for a contradiction. Put $w_{\l}(x)=\frac{v_{\l}(x)}{|v_{\l}|_{L^{\infty}}}$. By (\ref{homog}), remark that $w_{\l}$ satisfies 
$$
	\l w_{\l}+F(x,\n w_{\l},\n^2 w_{\l})
	 -\int_{{\bf R^N}}  
	[w_{\l}(x+z)-w_{\l}(x)
	\qquad\qquad\qquad
$$
\begin{equation}\label{wl}
	-{\bf 1}_{|z|\leq 1}\la \n w_{\l}(x),z \ra] c(z) dz 
	-\frac{f(x)-\l u_\l (x_0)}{|v_{\l}|_{L^{\infty}}}=0 \qquad x\in {\bf T^N},
\end{equation}
and that 
\begin{equation}\label{nonzero}
	|w_{\l}|_{L^{\infty}}=1, \quad w_{\l}(x_0)=0 \qquad \forall \l\in (0,1).
\end{equation}
From the comparison result for (\ref{ax}), there exists a constant $C>0$ such that $|\l u_{\l}|_{L^{\infty}}<C$ $(\forall \l\in (0,1))$, and thus from (\ref{increase}) and (\ref{nonzero}) there exists a constant $M>0$ such that 
$$
	|\l' w_{\l'}
	-\frac{f(x)-\l' u_{\l'} (x_0)}{|v_{\l'}|_{L^{\infty}}}|\leq M,
	\quad
	|w_{\l'}|\leq M. 
$$
Therefore, by applying Thorems 3.1 and 3.2 to (\ref{wl}) for $g=-(\l' w_{\l'}-\frac{f(x)-\l' u_{\l'} (x_0)}{|v_{\l'}|_{L^{\infty}}})$ and $\theta_0=0$, we know that there exist $\theta\in (0,1)$ and a  constant $C_{\theta}>0$ such that 
$$
	|w_{\l'}(x)-w_{\l'}(y)|\leq C_{\theta} |x-y|^{\theta} \quad 
	\forall x,y\in {\bf T^N}.
$$
So, by the Ascoli-Alzera theorem, 
there exists an H{\"o}lder continuous function 
$w$ such that 
$$
	\lim_{\l'\to 0} w_{\l'}(x)=w(x) \quad \hbox{uniformly in}\quad 
	{\bf T^N}, 
$$
and from (\ref{nonzero})
$$
	|w|_{L^{\infty}}=1,\quad w(x_0)=0.
$$
Moreover, by putting $\l=\l'$ in (\ref{wl}), and by passing $\l'\to 0$, 
since the limit procedure of viscosity solutions, introduced by Barles and Perthame \cite{bp1} 
(see also \cite{ar4} and \cite{users}), is valid for the present nonlocal case, we see that $w$ is a viscosity solution of 
$$
	F(x,\n w,\n^2 w)-\int_{{\bf R^N}}  
	[w(x+z)-w(x)
	\qquad\qquad\qquad\qquad\qquad\qquad
$$
\begin{equation}\label{wlimit}
	\qquad\qquad
	-{\bf 1}_{|z|\leq 1}\la \n w(x),z \ra] c(z) dz =0 \qquad \forall x\in {\bf T^N}.
\end{equation}
However, since $F$ satisfies (\ref{positive}), by the strong maximum principle 
 (Theorem 5.1), and by taking account that $w$ is periodic in ${\bf T^N}$, we see that 
 $w$ is almost everywhere constant in ${\bf T^N}$. This contradicts to the fact that $w$ is an H{\"o}lder continuous function such that $|w|_{L^{\infty}}=1$ and $w(x_0)=0$. 
 Therefore, the assumption (\ref{increase}) is false, and we have proved (\ref{M}). \\
 
(Step 2.) From Step 1, we see that $v_{\l}$ ($\l\in (0,1)$) satisfies (\ref{M})
 and 
$$	
	F(x,\n v_{\l},\n^2 v_{\l}) -\int_{{\bf R^N}}  
	[v_{\l}(x+z)-v_{\l}(x)
	\qquad\qquad\qquad\qquad\qquad
$$
\begin{equation}\label{vl}
	-{\bf 1}_{|z|\leq 1}\la \n v_{\l}(x),z \ra] c(z) dz 
	-(f-\l u_\l (x_0)-\l v_{\l})=0 \quad \hbox{in}\quad {\bf T^N},
\end{equation}
 and there exists $M>0$ such that
$$
	|f(x)-\l u_\l (x_0)-\l v_{\l}(x)|<M \quad \forall x\in {\bf T^N},
	\quad \forall \l\in (0,1).
$$
 We apply again the result in Theorems 3.1 and 3.2  to (\ref{vl}), and see that
 there 
 exist  $\theta\in (0,1)$ and a constant $C_{\theta}>0$ such that 
$$
	|v_{\l}(x)-v_{\l}(y)|\leq C_{\theta}|x-y|^{\theta}\quad \forall x,y\in {\bf T^N}, \quad \forall \l\in (0,1). 
$$
So, we can take a subsequence $\l'\to 0$ of $\l \to 0$ such that 
$$
	\lim_{\l'\to 0}v_{\l'}(x)=\lim_{\l'\to 0}(u_{\l'}(x)-u_{\l'}(x_0))
	= \exists u(x)\quad
	\hbox{uniformly in}\quad {\bf T^N},
$$
$$
	\lim_{\l'\to 0} \l' u_{\l'}(x)=\lim_{\l'\to 0} \l' u_{\l'}(x_0)=d_f
	\quad
	\hbox{uniformly in}\quad {\bf T^N}.
$$
In the next step, we shall prove that the limit $d_f$ is independent on the 
choise of the subsequence $\l'\to 0$.\\

(Step 3.) We shall prove the uniqueness of the limit number $d_f$ obtained in Step 2. Let $(d_f,u)$, and $(d'_f,u')$ $(d_f\neq d'_f)$ be 
two pairs of the limit numbers and the limit functions. Thus, 
$$	
	d_f+F(x,\n u,\n^2 u)-\int_{{\bf R^N}}  
	[u(x+z)-u(x)
	\qquad\qquad\qquad\qquad\qquad\quad
$$
$$
	\qquad\qquad
	-{\bf 1}_{|z|\leq 1}\la \n u(x),z \ra] c(z) dz
	-f(x)=0 \quad \hbox{in}\quad {\bf T^N}, 
$$
and 
$$	
	{d'}_f+F(x,\n u',\n^2 u')-\int_{{\bf R^N}}  
	[u'(x+z)-u'(x)
	\qquad\qquad\qquad\qquad\qquad
$$
$$
	\qquad\qquad
	-{\bf 1}_{|z|\leq 1}\la \n u'(x),z \ra] c(z) dz
	-f(x)=0 \quad \hbox{in}\quad {\bf T^N}. 
$$
We may assume that $d'_f< d_f$, and by adding a constant if necessary we may also assume that $u>u'$. 
For any small $\nu>0$, by choosing $\l>0$ small enough 
 we see that $u$ and $u'$ are respectively a viscosity subsolution and a viscosity supersolution of the following problems.
 $$	
	\l u +F(x,\n u,\n^2 u)-\int_{{\bf R^N}}  
	[u(x+z)-u(x)
	\qquad\qquad\qquad\qquad\qquad\qquad
$$
$$
	\qquad\qquad
	-{\bf 1}_{|z|\leq 1}\la \n u(x),z \ra] c(z) dz
	-f(x)\leq \nu-d_f \quad \hbox{in}\quad {\bf T^N}. 
$$
$$	
	\l u' +F(x,\n u',\n^2 u')-\int_{{\bf R^N}}  
	[u'(x+z)-u'(x)
	\qquad\qquad\qquad\qquad\qquad\qquad
$$
$$
	\qquad\qquad
	-{\bf 1}_{|z|\leq 1}\la \n u'(x),z \ra] c(z) dz
	-f(x)\geq -\nu-d'_f \quad \hbox{in}\quad {\bf T^N}. 
$$
Then, from the comparison result 
$$
	0<\l (u-u')(x) \leq d'_f-d_f<0,
$$
which is a contradiction. Thus, $d_f$ is the unique number such that (\ref{ergodic}) has a viscosity solution. We have proved the claim of (i).\\

(ii) We treat the case that the partial differential operator is $F$. (The proof for the case of $H$ is same, and we omit 
 it.) 
Let $v_{\l}=u_{\l}-u_{\l}(x_0)$, and put 
$
	|v_{\l}|_{\infty}=\frac{C_{\l}}{\l}.
$
We shall prove the claim in the following three steps.\\
(Step 1.) If for a subsequence $\l'\to 0$, $\lim_{\l'\to 0}C_{\l'}= 0$, then 
$$
	|\l' u_{\l'}(x)-\l' u_{\l'}(x_0)|_{\infty}=\l' |v_{\l'}|_{\infty}=C_{\l'}
	\to 0,
$$
which implies the existence of a constant 
$$
	d_f=\lim_{\l' \to 0}\l' u_{\l'}(x)=\lim_{\l' \to 0}\l' u_{\l'}(x_0)
	\qquad \hbox{uniformly in} \quad {\bf T^N}.
$$

(Step 2.) Now, assume that for any subsequence $\l'\to 0$, $C_{\l'}$ does not converge to zero. That is, there exists a number $C_0>0$ such that $\lim\inf_{\l\to 0}C_{\l}\geq C_0>0$. From the comparison result for  (\ref{ax}), $|v_{\l}|_{\infty}\leq \frac{2M}{\l}$, and thus $0<C_{\l}\leq 2M$ ($\forall \l\in (0,1)$) holds. Hence, we can take 
a subsequence 
 $\l'\to 0$ such that $\lim_{\l'\to 0} C_{\l'}=\overline{C}$ 
 ($\lim_{\l'\to 0}\l'|u_{\l'}-u_{\l'}(x_0)|_{\infty}=\overline{C}$), where $C_0\leq \overline{C}\leq 2M$. For simplicity, we shall use $\l$ in place of $\l'$. Then, for 
 $w_{\l}=\frac{v_{\l}}{|v_{\l}|_{\infty}}$ we have 
\begin{equation}\label{wlC}
	\l w_{\l} +F(x,\n w_{\l},\n^2 w_{\l})-\int_{{\bf R^N}}  [w_{\l}(x+z)-w_{\l}(x)
	\qquad\qquad\qquad\qquad
$$
$$
	-\la \n w_{\l}(x),z \ra] c(z) dz 
	-\frac{\l}{C_{\l}}(f(x)-\l u_{\l}(x_0))=0 \quad \hbox{in}\quad {\bf T^N},
\end{equation}
where $|w_{\l}|_{\infty}=1$, $w_{\l}(x_0)=0$, $\lim_{\l\to 0}C_{\l}=\overline{C}$. 
	Here, we claim that for a constant $\theta\in (0,\theta_0)$, there exists $C_{\theta}>0$ independent on $\l\in (0,1)$ such that 
\begin{equation}\label{whold}
	|w_{\l}(x)-w_{\l}(y)|\leq C_{\theta}|x-y|^{\theta},
\end{equation}
	which can be  proved by the similar contradiction argument used in the proof of Theorem 4.2, and which we omit here.

	From (\ref{whold}), there exists a subsequence $\l'\to 0$ such that 
$\lim_{\l'\to 0}w_{\l'}(x)=\exists w(x)$, where the limit $w$ is also H{\"o}lder continuous, $|w|_{\infty}=1$, $w(x_0)=0$ and is the viscosity solution of 
$$
	F(x,\n w,\n^2 w) -\int_{{\bf R^N}}  
	[w(x+z)
	-w(x)-\la \n w(x),z \ra] c(z) =0 \quad \hbox{in}\quad {\bf T^N}.
$$
However, the strong maximum principle (Theorem 5.1) asserts that $w$ is almost everywhere constant, which is a contradiction. Therefore, $\lim\inf_{\l'\to 0}C_{\l'}=\overline{C}>0$ is false. \\
(Step 3) From Steps 1 and 2, we see that there exists a subsequence $\l'\to 0$ such that $\lim_{\l'\to 0}\l' u_{\l'}(x)=d_f$ uniformly in ${\bf T^N}$. 
Therefore, for any $\nu>0$ there exists $\l'>0$ small enough such that 
$$
	d_f+F(x,\n u_{\l'},\n^2 u_{\l'}) -\int_{{\bf R^N}}  
	[u_{\l'}(x+z)
	\qquad\qquad\qquad\qquad\qquad\qquad\qquad
$$
$$
	-u_{\l'}(x)-{\bf 1}_{|z|\leq 1}\la \n u_{\l'}(x),z \ra] c(z)dz -f(x)\leq \nu \quad x\in {\bf T^N},
$$
\begin{equation}\label{apcell}
	\quad
\end{equation}
$$
	d_f+F(x,\n u_{\l'},\n^2 u_{\l'}) -\int_{{\bf R^N}}  
	[u_{\l'}(x+z)
	\qquad\qquad\qquad\qquad\qquad\qquad\qquad
$$
$$
	-u_{\l'}(x)-{\bf 1}_{|z|\leq 1}\la \n u_{\l'}(x),z \ra] c(z)dz -f(x)\geq -\nu \qquad x\in {\bf T^N}.
$$
The uniqueness of $d_f$ can be proved in a similar way to the proof for (i).\\

%%%%%%%%%%%%%%%%%%%%%%%%%%%%%%%%%%%%%%%%%%%%%%%%%%%%%%%%%%%%%%%%%%%%%%%%%%%

\end{document}